\documentclass[12pt,leqno,]{article}
\usepackage{amsmath,amsthm,amssymb}
\usepackage{color}
\usepackage[utf8]{inputenc}

\newtheorem{theorem}{THEOREM}
\newtheorem{lemma}{Lemma}
\newtheorem{remark}{REMARK}

\newcommand{\dx}{\, \mathrm{d}x}

\newcommand{\ds}{\, \mathrm{d}s}
\newcommand{\dr}{\, \mathrm{d}r}
\newcommand{\dt}{\, \mathrm{d}t}

\newcommand{\rz}{\mathbb{R}^{2}}

\newcommand{\R}{\mathbb{R}}
\renewcommand{\div}{\operatorname{div}}

\newcommand{\ep}{\varepsilon}

\def\IntQRx{\int_{{Q_R (x_0)}}}
\def\IntQzd{\int_{Q_{\frac{3}{2}R}(x_0)}}
\def\IntTzd{\int_{T_{\frac{3}{2}R}(x_0)}}
\def\IntT2R{\int_{T_{2R}(x_0)}}

\def\IntQ2Rx{\int_{{Q_{2R} (x_0)}}}

\def\Ints0{\int^s_0}
\def\Intt0{\int^t_0}

\newcommand{\loc}{\operatorname{loc}}

\renewcommand{\div}{\operatorname{div}}

\newcommand{\newmathop}[2]{\newcommand{#1}{\mathop{\mit{#2}}}}
\newmathop{\Minttext}{\int\limits\hspace{ -9.2em}-}
\newmathop{\Mint}{\int\hspace{ -1.0em}-}

\addtocounter{section}{0}

\title{\bf A Liouville theorem for stationary incompressible fluids of von Mises type}
\author{M.~Fuchs \and J.~M\"uller }
\date{}

\begin{document}

\parindent2ex

\maketitle
\begin{minipage}{16cm}
AMS Subject-Classification: 76D05, 76D07, 76M30, 35Q30 \\

\end{minipage}

\begin{abstract}
We consider entire solutions $u$ of the equations describing the stationary flow of a generalized Newtonian fluid in 2D concentrating on the question, if a Liouville-type result holds in the sense that the boundedness of $u$ implies its constancy. A positive answer is true for $p$-fluids in the case $p>1$ (including the classical Navier-Stokes system for the choice $p=2$), and recently we established this Liouville property for the Prandtl-Eyring fluid model, for which the dissipative potential has nearly linear growth. Here we finally discuss the case of perfectly plastic fluids whose flow is governed by a von Mises-type stress-strain relation formally corresponding to the case $p=1$. it turns out that, for dissipative potentials of linear growth, the condition of $\mu$-ellipticity with exponent $\mu<2$ is sufficient for proving the Liouville theorem.
\end{abstract} 

In this note we look at entire solutions $u: \rz \to \rz$ of the homogeneous equation
\begin{equation}
\label{G1}
- \div \left[T^D (\ep (u))\right] + u^k \partial_k u + \nabla \pi = 0
\end{equation}
together with the incompressibility condition
\begin{equation}
\label{G2}
\div u = 0 \, .
\end{equation}
Here $u$ denotes the velocity field of a fluid and $\pi : \rz \to \R$ is the a priori unknown pressure function. In equation (\ref{G1}) and in what follows we adopt the convention of summation with respect to indices repeated twice. By $\ep (u)$ we denote the symmetric gradient of the field $u$ and $T^D$ represents the deviatoric part of the stress tensor being characteristic for the fluid under consideration. For further mathematical and also physical explanations we refer to the monographs of 
Ladyzhenskaya \cite{La},  Galdi \cite{Ga1}, \cite{Ga2} and of M\'alek, Nec\v{a}s, Rokyta, R\r{u}\v{z}i\v{c}ka \cite{MNRR} as well as to the book \cite{FSe}. A case of particular interest arises, when $T^D$ is of the type
\begin{equation}
\label{G3}
T^D = \nabla H
\end {equation}
for a potential $H$ such that
\begin{equation}
\label{G4}
H (\ep) = h (|\ep|)
\end{equation}
holds with a given density $h : [0, \infty) \to [0, \infty)$ of class $C^2$. 
Combining (\ref{G3}) and (\ref{G4}) we see that 
\begin{equation}
\label{G5}
T^D (\ep) = \frac{h' (|\ep|)}{|\ep|} \ep
\end{equation}
holds, and equation (\ref{G5}) includes as particular cases

\begin{enumerate}
\item[(i)]  power law models:
\[
h (t) = t^p \ \mbox{or} \ = (\delta + t^2)^{p/2}, \ 1 < p < \infty, \ \delta > 0,
\]
\item[(ii)]
Prandtl-Eyring fluids: $h (t) = t \ln (1 + t), \ t \ge 0$.
\end{enumerate}
As a matter of fact we recover the Navier-Stokes equation (NSE) just by letting $h (t) = t^2$, and in their fundamental paper \cite{KNSS}, Koch, Nadirashvili, Seregin and Sver\v{a}k obtained the following Liouville-type result as a byproduct of their investigations on the regularity of solutions to the instationary variant of (NSE).
\begin{theorem}
Suppose that $u : \rz \to \rz$ is a solution of (\ref{G1}) and (\ref{G2}). Let in addition (\ref{G3}) and (\ref{G4}) hold for the choice $h (t) = \nu t^2$ with some positive constant $\nu$. Then, if
\begin{equation}
\label{G6}
\sup\limits_{\rz} |u| < \infty\, ,
\end{equation}
$u$ must be a constant vector.
\end{theorem}
For (NSE), different types of Liouville theorems were studied. For example, Gilbarg and Weinberger showed in their paper \cite{GW} the constancy of finite energy solutions to (NSE) in the plane, more precisely,  the conclusion of Theorem 1 remains valid if (\ref{G6}) is replaced by
\begin{equation}
\label{G7}
\int_{\rz} |\nabla u|^2 \dx < \infty \, .
\end{equation}
We wish to remark that the proofs of the above results use the linearity of the leading part of (NSE) in an essential way. However, referring to the results obtained in \cite{BFZ}, \cite{Fu}, \cite{FZ}, \cite{Zh1} and \cite{Zh2}, we could show by applying appropriate arguments:
\begin{theorem}
Suppose that $u \in C^1 (\rz, \rz)$ is a (weak) solution of equations (\ref{G1}) and (\ref{G2}) on the whole plane with $T^D$ given by (\ref{G3}) and (\ref{G4}) for a function $h$ being subject to (i) or (ii). Suppose that either (\ref{G6}) holds or that (\ref{G7}) is replaced by 
\begin{equation}
\label{G8}
\int_{\rz} h (|\nabla u|) \dx < \infty\, . 
\end{equation}
Then $u$ is a constant vector.
\end{theorem}
\begin{remark}
In the case of (NSE) weak solutions of some  local Sobolev class are automatically smooth. This is in general not clear for generalized Newtonian fluids and we therefore assume at least $u \in C^1$. A slightly weaker formulation can be found in \cite{JK}.
\end{remark}
\begin{remark}
A comprehensive survey of Theorem 2 including further related results is given in the paper \cite{Fu3}.
\end{remark}
\noindent Let us now turn to the von Mises flow (see, e.g., \cite{vMi}, \cite{HP}, \cite{Pr}) for which we formally have
\begin{equation}
\label{G9}
h (t) = t, \ t \ge 0\, ,
\end{equation}
which means (recall (\ref{G4})) $H (\ep (u)) = |\ep (u)|$. Since this density is neither differentiable nor strictly convex, equation (\ref{G3}) and thereby identity (\ref{G1}) can only be interpreted in a very weak sense, and we have no idea, if in this setting a Liouville-type result can be expected. For this reason we follow the ideas of \cite{BFM2} and replace the density from (\ref{G9}) through a family $h_\mu, \mu > 1$, of more regular densities being still of linear growth and such that $h_\mu (|\ep (u)|) \to |\ep (u)|$ as $\mu \to \infty$. For example we may take $(\mu > 1)$
\begin{equation}
\label{G10}
h_\mu (t) = (\mu - 1) \Phi_\mu (t), \ t \ge 0\, ,
\end{equation}
where we have abbreviated
\begin{equation}
\label{G11}
\begin{array}{l}
\hspace*{0.5cm}\Phi_\mu (t) :=\displaystyle{\Intt0} \displaystyle{\Ints0} (1 + r)^{- \mu} \dr \ds \\[3ex]
=\left\{
\begin{array}{l}
\frac{1}{\mu -1} t + \frac{1}{\mu - 1} \frac{1}{\mu - 2} (t + 1)^{-\mu + 2} - \frac{1}{\mu - 1} \frac{1}{\mu - 2}\, , \ \mu \not= 2\, ,\\[2ex]
t -\ln (1 + t), \ \mu = 2\, .
\end{array}
\right.
\end{array}
\end{equation}
Note that actually
\[
\lim_{\mu \to \infty} h_\mu (t) = t\, ,
\]
and if we formally let $\mu = 1$ in the first line of (\ref{G11}), then - up to negligible terms - we recover the Prandtl-Eyring model (ii). In the same spirit, the choice $\mu < 1$ leads to $p$-fluids with value $p : = 2 - \mu$. Of course our considerations are not limited to the particular density $h_\mu$. More general, we can choose any function $h : [0, \infty) \to [0, \infty)$ of class $C^2$ such that 
\begin{eqnarray}
\label{G12}
& \quad h (0) = h' (0) = 0\, ,\\[2ex]
\label{G13}
 & c_1 (1 + t)^{- \mu} \le \min \left\{h'' (t), \frac{h' (t)}{t}\right\}\, ,\\[2ex]
\label{G14}
& \quad \max \left\{ h'' (t), \frac{h' (t)}{t} \right\} \le \displaystyle\frac{c_2}{1 + t}
\end{eqnarray}
for any $t \ge 0$ and with  exponent
\begin{equation}
\label{G15}
\mu \in (1, \infty)\, ,
\end{equation}
$c_1, c_2$ denoting positive constants. It is immediate that the functions $h_\mu$ defined in (\ref{G10}) satisfy (\ref{G12}) - (\ref{G14}). Moreover, if we define $H$ according to (\ref{G4}), then the conditions (\ref{G12}) - (\ref{G14}) imply the $\mu$-ellipticity of $H$, i.e. 
\begin{equation}
\label{G16}
c_1 (1 + |\ep|)^{-\mu} |\sigma|^2 \le D^2 H (\ep) (\sigma, \sigma) \le c_2 (1 + |\ep|)^{-1}|\sigma|^2
\end{equation}
being valid for $(2 \times 2)$-matrices $\ep, \sigma$. In addition, the potential $H$ is of linear growth. More precisely we have (see Lemma 2.7 in \cite{BFM2})
\begin{equation}
\label{G17}
c_1 \ 2^{- \mu} (|\ep| - 1) \le H (\ep) \le c_2 |\ep|
\end{equation}
with $c_1$ and $c_2$ from (\ref{G13}) and (\ref{G14}), respectively. Now we can state our main result:
\begin{theorem}\label{Thm3}
Let $u \in C^1 (\rz, \rz)$ denote a (weak) solution of (\ref{G1}) and (\ref{G2}) on the whole plane with deviatoric stress tensor defined according to (\ref{G3}) and (\ref{G4}), where $h$ satisfies the conditions (\ref{G12}) - (\ref{G15}). In addition we assume that
\begin{equation}
\label{G18}
\mu < 2\, .
\end{equation}
If then the velocity field $u$ is bounded, it must be  a constant vector.
\end{theorem}
\begin{remark}
The limitation (\ref{G18}) enters for two reasons. First, as it will become evident from the proof of Theorem \ref{Thm3},  it plays the role of a technical restriction making it possible to manage certain quantities. Second, the results obtained in \cite{BFM2} suggest that there is some hope for the existence of regular weak solutions in case $\mu < 2$ motivating our assumption $u \in C^1$, whereas counterexamples taken from a slightly different setting (see \cite{FMT}) indicate that for $\mu > 2$ equations of the form (\ref{G1}) may fail to have solutions even on bounded domains, which can be found in some Sobolev space. Due to the linear growth of $H$ stated in (\ref{G17}) and with respect to the ellipticity condition (\ref{G16}), the space of functions with bounded deformation (see, e.g. \cite{ST2}, \cite{Su1}, \cite{Su2}) seems to be the appropriate class for discussing (\ref{G1}) but it is not evident how to give a reasonable ``very weak" formulation of equation (\ref{G1}) and to investigate the Liouville property in this setting.
\end{remark}

\begin{remark}
 We emphasize that our arguments for the proof of Theorem \ref{Thm3} fail in dimension $n\geq 3$. However, for the three dimensional Stokes-type problem,
\[
 \left\{\begin{aligned}-\div \left[\nabla H(\ep (u))\right] + \nabla \pi &= 0\\ \div u&=0\end{aligned}\right.
\]
we have a variant of Theorem 1.3 from \cite{FZ}, i.e. under the assumptions of Theorem \ref{Thm3} on the function $H$ and the parameter $\mu$,  every entire solution $u\in C^1(\R^3,\R^3)$ of the above system, for which $|x|^{-\alpha}|u(x)|$ is bounded with some $\alpha\in [0,1/2)$, must be constant. The proof follows along the lines of \cite{FZ}. In the case $n=2$, this result holds true even for the optimal parameter range $\alpha\in [0,1)$, cf. Theorem 1.1 in \cite{FZ}.

\end{remark}

\noindent {\bf Proof of Theorem \ref{Thm3}.} \ We follow the arguments outlined in Section 3 of \cite{FZ} keeping the notation introduced in Theorem \ref{Thm3} and assuming that all the hypothesis of Theorem \ref{Thm3} are valid. We start with
\begin{lemma}\label{Lem1}
There exists a constant $c = c \left(\|u\|_{L^\infty (\rz)}\right) < \infty$ such that
\begin{equation}
\label{G19}
\IntQRx H \left(\ep (u)\right) \dx \le c R
\end{equation}
holds for any square $Q_R (x_0) \subset \rz$, $Q_R (x_0) := \left\{ x \in \rz : |x_i - x_{0i}| < R, i = 1,2\right\}, \ R > 0, \ x_0 \in \rz$.
\end{lemma}
\noindent {\bf Proof of Lemma \ref{Lem1}.} \ From equations (\ref{G1}) and (\ref{G2}) we infer (recalling also (\ref{G3})) after an integration by parts
\begin{equation}
\label{G20}
0=\IntQ2Rx  D H \left(\ep (u)\right) : \ep (\varphi) \dx + \IntQ2Rx u^k \partial_k u^i \varphi^i \dx
\end{equation}
for any  test vector $\varphi$ such that $\div \varphi = 0$ on $Q_{2R} (x_0)$ and $\varphi = 0$ on $\partial Q_{2R} (x_0)$. Let $\eta\in C_0^\infty(Q_{2R}(x_0))$ be a cut-off function such that $0\leq \eta(x)\leq 1$, $\eta\equiv 1$ on $Q_R(x_0)$ and $|\nabla\eta|\leq cR^{-1}$, with $c$ denoting a generic positive constant. Observing that the scalar function $f=\div(\eta^2 u)$ fulfills the assumptions of Lemma 3.2 in \cite{FZ}, we infer that there exists $w\in \overset{\circ}{W}{}^{1,2}(Q_{2R}(x_0),\R^2)$ such that
\begin{align*}
&\div w=f\quad\text{together with}\\
&\|\nabla w\|_{L^2(Q_{2R}(x_0))}\leq c \|f\|_{L^2(Q_{2R}(x_0))}
\end{align*}
for a constant $c>0$ which is independent of $R$ and $x_0$. Now choosing $\varphi=\eta^2 u-w$ in (\ref{G20}), we arrive at the identity
\begin{align}\label{G21}
&\IntQ2Rx  D H \left(\ep (u)\right) : \ep (u) \eta^2 \dx + 2 \IntQ2Rx \frac{\partial H}{\partial \ep_{i \alpha}} \left(\ep (u)\right) \partial_\alpha \eta u^i \eta \dx\\
& - \IntQ2Rx DH \left(\ep (u)\right) : \ep (w) \dx  \nonumber \\
& =\IntQ2Rx u^k \partial_k u^i w^i \dx - \IntQ2Rx u^k \partial_k u^i u^i \eta^2 \dx.\nonumber
\end{align}
Let us look at the quantities on the l.h.s. of (\ref{G21}): for any tensor $\ep$ it holds (recall (\ref{G4}) and (\ref{G5}))
\[
D H (\ep) : \ep = \frac{h'(|\ep|)}{|\ep|} |\ep|^2 \, ,
\]
and the convexity of $h$ (see (\ref{G13})) together with (\ref{G12}) implies
\[
0 = h (0) \ge h (t) - t h' (t)\, ,
\]
thus $t h' (t) \ge h (t)$  and therefore
\begin{align}\label{G22}
 \IntQ2Rx DH \left(\ep (u)\right) : \ep (u) \eta^2 \dx \ge  \IntQ2Rx \eta^2 H \left(\ep (u)\right) \dx\, .
\end{align}
Using the boundedness of $D H$ (compare (\ref{G14})) we see
\[
\left| 2  \IntQ2Rx \frac{\partial H}{\partial \ep_{i \alpha}} \left(\ep (u)\right) \partial_\alpha \eta u^i \eta \dx \right|  \le  c  \IntQ2Rx |\nabla  \eta| |u| \dx \le c R\nonumber 
\]
on account of $|\nabla \eta| \le c/R$ by the choice of $\eta$ and due to the boundedness of $u$. Again by the boundedness of $D H$ it follows from the properties of $w$
\begin{align}
&\left| \IntQ2Rx D H \left(\ep (u)\right) : \ep (w) \dx \right| \le c \IntQ2Rx \left|\ep (w)\right| \dx  \nonumber\\[2ex]
& \le c R \left\| \ep (w)\right\|_{L^2 (Q_{2R} (x_0))} \le c R \left\| \div (\eta^2 u)\right\|_{L^2 (Q_{2 R} (x_0))} \nonumber\\[2ex]
& \le c R \left\| u\cdot \nabla \eta \right\|_{L^2 (Q_{2 R} (x_0))} \le c R\, .\nonumber
\end{align}
Let us note that during our calculations $c$ is a generic constant independent of $R$ and $x_0$. Combining the above estimates with (\ref{G22}) and returning to (\ref{G21}), we get
\begin{align}\label{G23}
\IntQ2Rx \eta^2 H \left(\ep (u)\right) \dx \le c R + \left| \mbox{r.h.s. of \ (\ref{G21})} \right|\, .
\end{align}
For the quantities occurring on the r.h.s. of (\ref{G21}) we can quote (3.6) and (3.7) in \cite{FZ}, hence
\begin{align}\label{G24}
\left|\mbox{r.h.s. of \ (\ref{G21})}\right| \le c R \, ,
\end{align}
and by inserting (\ref{G24}) into (\ref{G23}) the claim of Lemma \ref{Lem1} follows. \qed

\noindent Up to now neither the condition of $\mu$-ellipticity (see (\ref{G13}) and (\ref{G16})) nor the bound (\ref{G18}) on the parameter $\mu$ have entered our discussion, which means that estimate (\ref{G19}) is valid under much weaker hypotheses as required in Theorem \ref{Thm3}. The full strength of our assumptions is needed in the next step. We return to equation (\ref{G20}) and replace $\varphi$ by $\partial_\alpha \varphi$ for $\varphi \in C^\infty_0 (Q_{\frac{3}{2} R} (x_0), \rz)$, $\div \varphi = 0$. After an  integration by parts we obtain

\begin{align}\label{G25}
&0 = \int_{Q_{\frac{3}{2} R} (x_0)} D^2 H \left(\ep (u)\right) \left(\partial_\alpha \ep (u), \ep (\varphi)\right) \dx\\
&- \int_{Q_{\frac{3}{2} R} (x_0)} u^k \partial_k u^i \partial_\alpha \varphi^i \dx, \quad \alpha = 1, 2\, . \nonumber
\end{align}
At this stage we recall our assumption $u \in C^1 (\rz, \rz)$, which enables us with the help of the difference-quotient technique to deduce $u \in W^{2,2}_{\loc} (\rz, \rz)$ from equation (\ref{G20}) and to justify (\ref{G25}). Note also that (\ref{G25}) corresponds exactly to formula (3.10) in \cite{FZ} and with the choice of $\varphi$ as in this reference, we deduce the identity (3.12). Without any changes we can pass to inequality (3.18) from this reference, i.e. setting $\omega := D^2 H \left(\ep (u)\right) \left(\partial_\alpha \ep (u), \partial_\alpha \ep (u) \right)$, we have
\begin{align}\label{G25'}
 \intop_{Q_{\frac{3}{2}R}(x_0)}\varphi^2\omega\dx\leq \delta \intop_{Q_{\frac{3}{2}R(x_0)}}\omega\dx+ c(\delta)\left[\frac{1}{R^2}\intop_{Q_{\frac{3}{2}R}(x_0)}|\nabla u|^2\dx+\frac{1}{R}\intop_{Q_{\frac{3}{2}R}(x_0)}|\nabla u|^2\dx\right]
\end{align}
for any  $\delta > 0$ and all squares $Q_{2R} (x_0)$. Specifying $\varphi$ as $0 \le \varphi \le 1$, $\varphi = 1$ on $Q_{\frac{3}{2} R} (x_0)$ and $|\nabla \varphi| \le c/R$, we my furthermore pass to inequality (3.19) from \cite{FZ} which reads as
\begin{align}\label{G26}
&\IntQRx \omega \dx \le \delta \IntQ2Rx \omega \dx \\
& + c (\delta) \left[R^{-4} \IntQ2Rx |u|^2 \dx + R^{-3} \IntQ2Rx |u|^2 \dx\nonumber \right.\\
& \left.+ R^{-2} \IntQ2Rx  \varphi^2 |\ep (u)|^2 \dx + R^{-1} \IntQ2Rx \varphi^2 |\ep (u)|^2 \dx \right].\nonumber
\end{align}
 We see that it remains to discuss the quantity $\IntQ2Rx \varphi^2 |\ep (u)|^2 \dx $. We have
\begin{align}
&\IntQ2Rx \varphi^2 \left|\ep (u)\right|^2 \dx = - \IntQ2Rx u^i \partial_j \left(\ep_{ij} (u)\varphi^2\right) \dx \nonumber\\
&=- \IntQ2Rx u^i  \partial_j \ep_{ij} (u) \varphi^2 \dx - \IntQ2Rx u^i \ep_{ij} (u) \partial_j \varphi^2 \dx \nonumber\\
&\le c \left[\IntQ2Rx \left|\nabla \ep (u)\right| \varphi^2 \dx + R^{-1} \IntQ2Rx  \left|\ep (u)\right|\dx \right]\nonumber
\end{align}
and (for any $\tau > 0$)
\begin{align}\label{G26'}
& \IntQ2Rx \left|\nabla \ep (u)\right| \varphi^2 \dx =  \IntQ2Rx  \varphi \left|\nabla \ep (u) \right| \left(1 + |\ep (u)|\right)^{-\mu/2} \left(1 + |\ep (u)|\right)^{\mu/2} \varphi \dx\\
&\le \tau \IntQ2Rx \omega \dx + c \tau^{-1} \IntQ2Rx \varphi^2 \left(1 + |\ep (u)|\right)^\mu \dx \,\nonumber , 
\end{align}
where  we have used H\"older's inequality and condition (\ref{G16}). Let us choose $\tau := \delta c (\delta)^{-1} R^2$ in (\ref{G26'}) with $c (\delta)$ from (\ref{G26}). Then it follows with a new constant $\tilde{c} (\delta)$
\begin{align}\label{G27}
&c (\delta) R^{-2} \IntQ2Rx \varphi^2 \left|\ep (u)\right|^2 \dx \le \delta \IntQ2Rx \omega \dx +\\
&\tilde{c} (\delta) \left[R^{-4} \IntQ2Rx \varphi^2 \left(1 + \left|\ep (u)\right|\right)^{\mu} \dx + R^{-3} \IntQ2Rx \left|\ep (u)\right|\dx \right]\, , \nonumber
\end{align}
and if we select $\tau := \delta c (\delta)^{-1} R$ in (\ref{G26'}), we find
\begin{align}\label{G28}
&c (\delta) R^{-1} \IntQ2Rx \varphi^2 \left|\ep (u)\right|^2 \dx \le \delta \IntQ2Rx \omega \dx +\\
&\tilde{c} (\delta) \left[R^{-2} \IntQ2Rx \varphi^2 \left(1 + \left|\ep (u)\right|\right)^{\mu} \dx + R^{-2} \IntQ2Rx \left|\ep (u)\right|\dx \right]\, . \nonumber
\end{align}
In (\ref{G27}) and (\ref{G28}) we have to get rid of the quantity $\IntQ2Rx  \varphi^2 (1 + |\ep (u)|)^\mu \dx$ by absorbing it into the left-hand sides.
\noindent It holds for any $\lambda > 0$
\begin{align}\label{G27'}
&\int_{Q_{2R}(x_0)} \left( 1 + |\ep (u)|\right)^{\mu} \varphi^2 \dx \le c \left[R^2+\IntQ2Rx \left|\ep (u)
\right|^\mu \varphi^2 \dx \right]  \\
&\leq c \left[ R^2 + \lambda \IntQ2Rx \varphi^2 |\ep (u)|^2 \dx + R^2 \lambda^{\frac{\mu}{\mu - 2}} \right]\nonumber
\end{align}
on account of Young's inequality and due to our assumption $\mu < 2$. By choosing $\lambda$ proportional to $R^2$ we infer from (\ref{G27})
\begin{align}\label{G29}
&c (\delta) R^{-2} \IntQ2Rx \varphi^2 \left|\ep (u)\right|^2 \dx \le \delta \IntQ2Rx \omega \dx \\
&+\tilde{c} (\delta) \left[R^{-2} + R^{{-2} + \frac{2\mu}{\mu - 2}} + R^{- 3} \IntQ2Rx |\ep (u)|\dx \right]\, ,\nonumber
\end{align}
whereas $\lambda \sim R$ in combination with (\ref{G28}) yields 
\begin{align}\label{G30}
&c (\delta) R^{-1} \IntQ2Rx \varphi^2 |\ep (u)|^2 \dx \le \delta \IntQ2Rx \omega \dx + \\
&\tilde{c} (\delta) \left[1 + R^{\frac{\mu}{\mu - 2}} + R^{-2} \IntQ2Rx |\ep (u)| \dx \right] \, .\nonumber
\end{align}
We insert (\ref{G29}) and (\ref{G30}) into (\ref{G26}), replacing $\delta$ by $\delta/3$ and get for any $\delta > 0$ and all squares $Q_{2R} (x_0) \subset \rz$
\begin{align}\label{G31}
&\IntQRx \omega \dx \le \delta \IntQ2Rx \omega \dx  \\
&+c (\delta) \left[ R^{-4} \IntQ2Rx |u|^2 \dx +R^{-3} \IntQ2Rx |u|^2 \dx  \right.\nonumber\\
& \quad \quad \quad +R^{-4} \IntQ2Rx 1 \dx + R^{- 4 + \frac{2\mu}{\mu - 2}} \IntQ2Rx 1 \dx \nonumber\\
&\quad \quad  \quad +R^{-2} \IntQ2Rx 1 \dx + R^{- 2 + \frac{\mu}{\mu - 2}} \IntQ2Rx 1 \dx  \nonumber\\
& \quad \quad \quad +\left. R^{-3} \IntQ2Rx \left|\ep (u)\right| \dx + R^{-2} \int_{Q_{2R}(x_0)} \left| \ep (u)\right| \dx \right]\,  \nonumber\\
&\quad \quad \quad=:\delta \IntQ2Rx \omega\dx+c(\delta)\Theta(R)\nonumber.
\end{align}
To inequality (\ref{G31}) we can apply Lemma 3.1 from \cite{FZ} to get
\begin{align}\label{G32}
\int_{Q_{R}(x_0)} \omega \dx \le c\Theta(R)
\end{align}
for arbitrary squares $Q_R (x_0)$. Clearly, $\Theta(R)$ is bounded and hence
\begin{align}\label{omegainf}
\intop_{\R^2} \omega \dx=: \omega_\infty<\infty.
\end{align}
Our goal is to show $\Theta(R)\to 0$ as $R\to\infty$ since, together with inequality (\ref{G32}), this implies $\ep(\nabla u)\equiv 0$, hence $\nabla^2u\equiv 0$ which means that $u$ is an affine function. The assumption on the boundedness of $u$ then gives the assertion of Theorem \ref{Thm3}.

We start with the term $R^{-2}\int_{Q_{2R}(x_0)}|\ep(u)|\dx$ by noting that due to $h'(0)=h(0)=0$ we have (for $\beta>0$ arbitrarily small)
\begin{align*}
h(x)&=\intop_0^1(1-t)h''(tx)\dt x^2\overset{(\ref{G13})}{\geq}c_1x^2\intop_0^1(1-t)(1+tx)^{-\mu}\dt\geq c_1x^2\intop_0^{\beta/x}\frac{1-t}{(1+tx)^\mu}\dt\\
&\geq c_1x^2\intop_0^{\beta/x}(1-t)\dt\,(1+\beta)^{-\mu}=c_1x^2(1+\beta)^{-\mu}\Big(\frac{\beta}{x}-\underset{\leq \frac{1}{2}\beta\frac{1}{x}}{\underbrace{\frac{1}{2}\frac{\beta^2}{x^2}}}\Big)\geq \frac{1}{2}\beta x c_1(1+\beta)^{-\mu},
\end{align*}
which means that
\begin{align}\label{G34}
h(x)\geq c(\beta)x\quad\text{for all }x\geq \beta.
\end{align}
Consequently, it holds
\[
R^{-2}\IntQ2Rx|\ep(u)|\dx=R^{-2}\intop_{\substack{Q_{2R}(x_0)\\ \cap\{|\ep(u)|\geq \beta\}}}|\ep(u)|\dx+R^{-2}\intop_{\substack{Q_{2R}(x_0)\\ \cap\{|\ep(u)|< \beta\}}}|\ep(u)|\dx=:T_1+T_2,
\]
with $\lim\limits_{R\to\infty}T_1=0$ by (\ref{G34}) and (\ref{G19}). Moreover, $\limsup\limits_{R\to\infty}T_2\leq c\beta$ and since $\beta$ may be chosen arbitrarily small, we infer 
\begin{align}\label{epstonull}
\lim\limits_{R\to\infty}\frac{1}{R^2}\IntQ2Rx|\ep(u)|\dx=0.
\end{align}

Next we note that, choosing a test function $\varphi\in C_0^\infty(Q_{2R}(x_0))$ with $\varphi\equiv 1$ on $T_{\frac{3}{2}R}(x_0):=Q_{\frac{3}{2}R}(x_0)-\overline{Q_R(x_0)}$ and such that $\operatorname{spt}\varphi\subset T_{2R}(x_0):=Q_{2R}(x_0)-\overline{Q_{\frac{R}{2}}(x_0)}$, we can derive the following more refined version of inequality (\ref{G25'}), which is identical to inequality (3.27) from \cite{FZ}:
\begin{align}\label{3.27'}
\intop_{Q_{R}(x_0)} \omega\dx\leq \delta \intop_{Q_{\frac{3}{2}R}(x_0)}\omega\dx&+c(\delta)\left[\frac{1}{R^2}\IntTzd |\nabla u|	^2\dx+\frac{1}{R}\IntTzd |\nabla u|^2\dx\right.\\
&\left.+\frac{1}{R}\left(\IntQzd |\nabla u|^2\dx\right)^{1/2}\left(\IntTzd |\nabla u|^2\dx\right)^{1/2}\right]. \nonumber
\end{align}
With the same arguments as in \cite{FZ}, we derive inequalites (3.28), (3.29) for the term $\IntTzd |\nabla u|^2\dx$:
\begin{align}\label{3.29'}
\begin{split}
\IntTzd |\nabla u|^2\dx&\leq c \left[\IntT2R\varphi^2|\ep(u)|^2\dx+\frac{1}{R^2}\int_{T_{2R}(x_0)}|u|^2\dx\right]\\
&\leq c\left[\Phi(R)+\frac{1}{R^2}\int_{T_{2R}(x_0)}|u|^2\dx\right],
\end{split}
\end{align}
where
\[
\Phi(R):=\left(\int_{T_{2R}(x_0)}\omega\dx\right)^{1/2}\left(\int_{T_{2R}(x_0)}\varphi^2\big(1+|\ep(u)|\big)^\mu\dx\right)^{1/2}+\frac{1}{R}\int_{T_{2R}(x_0)}|\ep(u)|\dx.
\]
With the estimate (\ref{G27'}) we infer
\begin{align*}
\Phi(R)\leq c \left(\int_{T_{2R}(x_0)}\omega\dx\right)^{1/2}\left[R+\left(\int_{T_{2R}(x_0)}\varphi^2|\ep(u)|^2 \dx\right)^{1/2}\right]+\frac{1}{R}\int_{T_{2R}(x_0)}|\ep(u)|\dx,
\end{align*}
and an application of Young's inequality  yields ($\alpha>0$ may be arbitrarily small)
\begin{align*}
\Phi(R)\leq  c \left(\int_{T_{2R}(x_0)}\omega\dx\right)^{1/2}\left[R+\alpha^{-1}+\alpha\int_{T_{2R}(x_0)}\varphi^2|\ep(u)|^2 \dx\right]+\frac{1}{R}\int_{T_{2R}(x_0)}|\ep(u)|\dx.
\end{align*}
Since 
\begin{align}\label{tonull}
\int_{T_{2R}(x_0)}\omega\dx\to 0\quad\text{ for }R\to\infty
\end{align}
 by (\ref{omegainf}), we may absorb the term $\alpha\int_{T_{2R}(x_0)}\varphi^2|\ep(u)|^2 \dx$ in the middle-term of the estimate (\ref{3.29'}) for $\alpha$ sufficiently small, thus
\begin{align*}
&\IntTzd |\nabla u|^2\dx\\
&\leq c \left[\left(\int_{T_{2R}(x_0)}\omega\dx\right)^{1/2}\left(R+\alpha^{-1}\right)+\frac{1}{R}\int_{T_{2R}(x_0)}|\ep(u)|\dx+\frac{1}{R^2}\int_{T_{2R}(x_0)}|u|^2\dx\right]
\end{align*}
and the boundedness of $u$ together with (\ref{tonull}) and (\ref{epstonull}) therefore implies
\[
\frac{1}{R}\IntTzd |\nabla u|^2\dx\to 0\quad\text{ for }R\to\infty.
\] 
It remains to discuss the quantity 
\[
\Psi(R):=\left(\frac{1}{R}\IntQzd |\nabla u|^2\dx\right)^{1/2}\left(\frac{1}{R}\IntTzd |\nabla u|^2\dx\right)^{1/2}.
\]
We have already established that the second factor goes to $0$ for $R\to\infty$ and it therefore suffices to show that the first factor is bounded. Arguing as above, we see that
\begin{align*}
&\IntQzd |\nabla u|^2\dx\\
&\leq c\left[\left(\IntQ2Rx\omega\dx\right)^{1/2}\left(R+\alpha^{-1}\right)+\frac{1}{R}\IntQ2Rx|\ep(u)|\dx+\frac{1}{R^2}\IntQ2Rx|u|^2\dx\right].
\end{align*}
Now, the boundedness of $\IntQ2Rx\omega\dx$ together with (\ref{epstonull}) clearly implies the boundedness of $R^{-1}\IntQzd |\nabla u|^2\dx$, which finishes the proof of Theorem \ref{Thm3}. \qed


\vspace*{3cm}
\vspace{3ex} 

\begin{tabular}{cc}
\begin{minipage}{7cm}
Jan M\"uller\\
Universit\"at des Saarlandes\\
Fachbereich 6.1 Mathematik\\
Postfach 15 11 50\\
D--66041 Saarbr\"ucken,\\ Germany\\
e--mail: jmueller@math.uni-sb.de
\end{minipage}&
\begin{minipage}{7cm}
Martin Fuchs\\
Universit\"at des Saarlandes\\
Fachbereich 6.1 Mathematik\\
Postfach 15 11 50\\
D--66041 Saarbr\"ucken\\ Germany\\
e--mail: fuchs@math.uni-sb.de
\end{minipage} 
\end{tabular}

\end{document}